\documentclass[11pt,reqno]{amsart}
\usepackage{amssymb,latexsym}
\usepackage{url}		

\usepackage{hyperref}

\calclayout

\usepackage[utf8]{inputenc}
\usepackage[T1]{fontenc}
\usepackage{amsfonts}
\usepackage[greek,english]{babel}
\usepackage{lmodern}

\usepackage{enumerate}
\usepackage{enumitem}
\usepackage{comment}

\usepackage[normalem]{ulem} 

\setbox0=\hbox{$+$}
\newdimen\plusheight
\plusheight=\ht0
\def\+{\mathbin{\lower\plusheight\hbox{$+$}}}

\setbox0=\hbox{$-$}
\newdimen\minusheight
\minusheight=\ht0
\def\cminus{\mathbin{\lower\plusheight\hbox{$-$}}}

\setbox0=\hbox{$\cdots$}
\newdimen\cdotsheight
\cdotsheight=\plusheight
\def\cds{\lower\cdotsheight\hbox{$\cdots$}}

\setbox0=\hbox{$\cdots$}
\newdimen\cdotsheight
\cdotsheight=\plusheight
\def\dds{\lower\cdotsheight\hbox{$\ddots$}}

\newcommand{\df}{\dfrac}

\theoremstyle{plain}
\numberwithin{equation}{section}
\newtheorem{thm}{Theorem}[section]
\newtheorem{theorem}[thm]{Theorem}

\makeatletter
\newcommand{\leqnomode}{\tagsleft@true\let\veqno\@@leqno}
\newcommand{\reqnomode}{\tagsleft@false\let\veqno\@@eqno}
\makeatother

\begin{document}
\setcounter{page}{1}

\title{The Rogers--Ramanujan continued fraction}
\author{Bruce C. Berndt and \"{O}rs Reb\'{a}k}
\address{Department of Mathematics, University of Illinois, 1409 West Green Street, Urbana, IL 61801, USA}
\email{berndt@illinois.edu}
\address{Department of Mathematics and Statistics, University of Troms\o{} -- The Arctic University of Norway, 9037 Troms\o{}, Norway}
\email{ors.rebak@uit.no}

\begin{abstract}
The primary purpose of this paper is to provide a survey of properties, values, identities, and generalizations of the Rogers--Ramanujan continued fraction, which is closely related to the Rogers--Ramanujan identities. Many of these results are found in Ramanujan's first two letters to Hardy, Ramanujan's notebooks, and his lost notebook. Short historical accounts are provided for both the notebooks and lost notebook.
\end{abstract}

\keywords{Rogers--Ramanujan continued fraction, theta functions, Ramanujan's lost notebook}
\subjclass[2020]{11--02, 11A55, 11F27, 05A30, 11F32}

\maketitle

\section{Introduction to Continued Fractions}

A (finite) continued fraction is an expression of the form
\begin{equation*}
b_0 + \cfrac{a_1}{b_1 + \cfrac{a_2}{b_2 + \cfrac{a_3}{b_3 + \ddots + \cfrac{a_n}{b_n}}}}
= b_0 + \df{a_1}{b_1} \+ \df{a_2}{b_2} \+ \df{a_3}{b_3} \+ \cds \+ \df{a_{n}}{b_n}.
\end{equation*}
For example,
\begin{equation}\label{ex1}
1 + \cfrac{1}{1 + \cfrac{2}{1 + \cfrac{3}{1}}}
 = 1+\df{1}{1}\+\df{2}{1}\+\df{3}{1} = \df{5}{3}.
\end{equation}
Note that on the right-hand sides above, the plus signs + appear, so to speak, in the `denominators' to avoid confusion with sums, where the + signs appear midway between successive fractions.

Suppose that we write
\begin{equation}\label{cf1}
b_0 + \cfrac{a_1}{b_1 + \cfrac{a_2}{b_2 + \cfrac{a_3}{b_3 + \ddots}}}
= b_0 + \df{a_1}{b_1} \+ \df{a_2}{b_2} \+ \df{a_3}{b_3} \+ \cds.
\end{equation}
The dots $\cdots$ indicate that the process does not stop, i.e., \eqref{cf1} is an \emph{infinite continued fraction}. As an example, let $\phi := (\sqrt{5}+1)/2$, the golden ratio. Observe that $\phi = 1 + 1/\phi$. Expanding recursively, we obtain the continued fraction
\begin{equation}\label{golden}
\phi = \df{\sqrt{5}+1}{2} = 1+ \df{1}{1} \+ \df{1}{1} \+ \df{1}{1} \+ \cds.
\end{equation}
Infinite continued fractions are considerably more interesting than finite continued fractions. And, they are even more interesting when the `numerators' and the `denominators' of the continued fractions are functions of a real or complex variable. In particular, in analogy with power series, the `numerators' of a continued fraction might be successive powers of $x$. To determine when such expressions \eqref{cf1} make sense, we need to establish theorems on the convergence and divergence of continued fractions, in analogy with theorems on the convergence and divergence of infinite series. For theorems on the convergence of continued fractions, consult the excellent book by Lisa Lorentzen and Haakon Waadeland \cite{lisa}.

Returning to \eqref{ex1}, we remark that the elegant infinite continued fraction
\begin{equation}\label{cf2}
\df{1}{1}\+\df{1}{1}\+\df{2}{1}\+\df{3}{1}\+\df{4}{1}\+\cds
\end{equation}
does converge, in contrast to the sum that is obtained by moving the plus signs `upward.' We shall discuss \eqref{cf2} at the conclusion of our paper.

\section{The Rogers--Ramanujan Continued Fraction}

One of the most famous continued fractions is the Rogers--Ramanujan continued fraction, defined by
\begin{equation}\label{rrcf}
R(q) := \df{q^{1/5}}{1} \+ \df{q}{1} \+ \df{q^2}{1} \+ \df{q^3}{1} \+ \cds, 
\qquad |q|<1,
\end{equation}
which was first studied by L.~J.~Rogers in an apparently unnoticed paper \cite{rogers-1894} at that time in 1894. Its alternating counterpart $S(q)$ is defined by \cite[p.~57]{geabcbI}
\begin{equation}\label{S}
S(q) := -R(-q) = \df{q^{1/5}}{1} \cminus \df{q}{1} \+ \df{q^2}{1} \cminus \df{q^3}{1} \+ \cds, \qquad |q|<1.
\end{equation}
Here, and in the sequel, we take the principal branch of $z^a$, for any complex number $z$ and rational number $a$. Readers might wonder why the factor $q^{1/5}$ is present on the right side of \eqref{rrcf}. Wouldn't the definition of $R(q)$ be simpler without it? With the factor $q^{1/5}$, $R(q)$ is a modular form, and consequently the theory of modular forms can be invoked to study $R(q)$.

The continued fraction $R(q)$ converges for all $q$ on the interior of the unit disc, i.e., for $|q|<1$. We shall show that the definition can be extended for some values with $|q| = 1$. Recalling \eqref{golden}, we see that \eqref{rrcf} converges at $q = 1$ to the reciprocal of the golden ratio. Similarly, we find that at $q = 1$ the continued fraction \eqref{S} converges to the golden ratio $\phi$. Hence, by elementary means, we have shown that, respectively,
\begin{equation}\label{golden1}
R(1) = \frac{1}{\phi} = \frac{\sqrt{5} - 1}{2} \qquad \text{and} \qquad
S(1) = \phi = \df{\sqrt{5} + 1}{2},
\end{equation}
where we have taken the principal value $(-1)^{1/5} = -1$ in \eqref{S}. Thus, the condition $|q| < 1$ is not a necessary condition. However, for $|q| > 1$, the continued fraction \eqref{rrcf} diverges.

Does $R(q)$ converge for other points on the unit circle? We cannot completely answer this question. However, I.~Schur \cite{schur} proved the following theorem.
\begin{theorem}
Let $q$ denote a primitive $n$th root of unity. If $n$ is a multiple of 5, then $R(q)$ diverges. Suppose that $n$ is not a multiple of 5. Let $\lambda= \big(\frac{n}{5}\big)$, the Legendre symbol.
Furthermore, let $\rho$ denote the least positive residue of $n$ modulo 5. Then, if $n$ is not a multiple of 5,
\begin{equation*}
R(q) = \lambda q^{(\lambda\rho n-1)/5}R(\lambda).
\end{equation*}
\end{theorem}
Observe that, by \eqref{golden1}, $R(\lambda)$ exists. On page~383 in his first notebook \cite{nb}, Ramanujan asserted that $R(q)$ diverges when $n$ is a multiple of 5, but his evaluation when $n$ is not a multiple of 5 is incorrect. See \cite[pp.~35--36]{V} for more details.

In his first letter to G.~H.~Hardy on 16 January 1913, among nearly 70 theorems, identities, and conjectures, Srinivasa Ramanujan \cite[p.~xxvii]{RamanujanCollected}, \cite[p.~29]{bcbrar} offered the continued fraction \eqref{rrcf} at the value $e^{-2\pi}$, more precisely,
\begin{align}
R(e^{-2\pi})
&= \df{e^{-2\pi/5}}{1}\+\df{e^{-2\pi}}{1}\+\df{e^{-4\pi}}{1}\+\df{e^{-6\pi}}{1}\+\cds\notag\\
&= \sqrt{\df{5 + \sqrt{5}}{2}} - \df{\sqrt{5} + 1}{2}= 5^{1/4}\sqrt{\phi} - \phi\label{2},
\end{align}
wherein the golden ratio \eqref{golden} appears \emph{twice} in this remarkable identity. Ramanujan \cite[p.~xxvii]{RamanujanCollected}, \cite[p.~29]{bcbrar} also provided Hardy with the evaluation of $R(q)$ at $-e^{-\pi}$, namely, with the use of \eqref{S},
\begin{equation}\label{3}
S(e^{-\pi}) = \df{e^{-\pi/5}}{1} \cminus \df{e^{-\pi}}{1}\+\df{e^{-2\pi}}{1} \cminus \df{e^{-3\pi}}{1}\+\cds = \sqrt{\df{5 - \sqrt{5}}{2}} - \df{\sqrt{5} - 1}{2} = \frac{5^{1/4}}{\sqrt{\phi}} - \frac{1}{\phi}.
\end{equation}
(There is a misprint in \eqref{3} in \cite[p.~xxvii]{RamanujanCollected}; replace `3' by `2' in the first denominator on the right-hand side.) In his reply of 8 February 1913, Hardy \cite[pp.~47--48]{bcbrar} lists several ``noteworthy theorems,'' one of which is \eqref{2}.

Ramanujan further asserted that \cite[p.~xvii]{RamanujanCollected}, \cite[p.~29]{bcbrar}
\begin{equation}\label{4}
\df{1}{1}\+\df{e^{-\pi\sqrt{n}}}{1}\+\df{e^{-2\pi\sqrt{n}}}{1}\+\df{e^{-3\pi\sqrt{n}}}{1}\+\cds
\end{equation}
``can be exactly found if $n$ be any positive rational quantity.'' Later, in this article, we show that Ramanujan's claim is correct, but for a given value of $n$, the evaluation may be extraordinarily difficult to find. As we shall relate in Section~\ref{values}, to determine \eqref{4} for a specific value of $n$, it is necessary to ascertain the values of relevant \emph{class invariants}. These are algebraic numbers, and they are frequently difficult to determine.

In his second letter to Hardy, dated 27 February 1913, Ramanujan \cite[p.~57]{bcbrar} repeats his claim in \eqref{4}. He then explicitly records the case, $n = 20$, to wit,
\begin{equation}\label{5}
\df{1}{1}\+\df{e^{-2\pi\sqrt{5}}}{1}\+\df{e^{-4\pi\sqrt{5}}}{1}\+\df{e^{-6\pi\sqrt{5}}}{1}\+\cds
= e^{2\pi/\sqrt{5}}\left\{\df{\sqrt{5}}{1+\sqrt[5]{5^{3/4}\left(\frac{\sqrt{5}-1}{2}\right)^{5/2} - 1}}-\df{\sqrt{5} + 1}{2}\right\}.
\end{equation}

Prior to stating \eqref{3}, Ramanujan \cite[p.~xviii]{RamanujanCollected}, \cite[p.~57]{bcbrar} offered to Hardy a gorgeous identity for $R(q)$ involving the golden ratio. We quote Ramanujan as he recorded it.

``If
\begin{equation*}
F(x) := \df{1}{1}\+\df{x}{1}\+\df{x^2}{1}\+\df{x^3}{1}\+\df{x^4}{1}\+\df{x^5}{1}\+\cds,
\end{equation*}
then
\begin{equation}\label{symmetry}
\left\{\df{\sqrt{5} + 1}{2} + e^{-2\alpha/5}F(e^{-2\alpha})\right\}\left\{\df{\sqrt{5} + 1}{2} + e^{-2\beta/5}F(e^{-2\beta})\right\} = 
\df{5 + \sqrt{5}}{2}
\end{equation}
with the condition $\alpha\beta = \pi^2$.''

Observe the beautiful symmetry in the identity. Ramanujan loved symmetry, and he proved many identities under a hypothesis similar to $\alpha\beta = \pi^2$ above. The identity is also called a `modular relation,' because it is reminiscent of the transformation $z \mapsto -1/z$, $z \in \mathbb{C}$ and $\operatorname{Im}(z) > 0$, in the theory of modular forms. If we set $\alpha=\beta=\pi$ in \eqref{symmetry}, we obtain Ramanujan's evaluation \eqref{2}. The first published proof of \eqref{symmetry} was by G.~N.~Watson \cite{watson1929} in 1929.

In his reply on 26 March 1913 (coincidentally, the day on which Paul Erd\H{o}s was born) to Ramanujan's second letter, Hardy \cite[pp.~77--78]{bcbrar} concluded with an appeal:

``What I should like above all is a definite proof of some of your results concerning continued fractions of the type
\begin{equation*}
\df{x}{1}\+\df{x^2}{1}\+\df{x^3}{1}\+\cds;
\end{equation*}
and I am quite sure that the wisest thing you can do, in your own interests, is to let me have one as soon as possible.''

\section{The Rogers--Ramanujan Identities}

The Rogers--Ramanujan continued fraction is associated with the Rogers--Ramanujan identities. To state them, we first need to define
\begin{equation*}
(a;q)_0 := 1, \qquad (a;q)_n := \prod_{k=0}^{n-1}(1-aq^k), \qquad n \geq 1,
\end{equation*}
and, if $|q| < 1$,
\begin{equation*}
(a;q)_{\infty} := \lim_{n \to \infty}(a;q)_n.
\end{equation*}
The Rogers--Ramanujan functions $G(q)$ and $H(q)$, and the Rogers--Ramanujan identities are given by
\begin{align}
G(q) := \sum^{\infty}_{n=0}\dfrac{q^{n^{2}}}{(q;q)_{n}} &=
\dfrac{1}{(q;q^5)_{\infty}(q^4;q^5)_{\infty}}, \qquad |q| < 1,\label{rr1}\\
H(q) := \sum^{\infty}_{n=0}\dfrac{q^{n(n+1)}}{(q;q)_{n}} &=
\dfrac{1}{(q^2;q^5)_{\infty}(q^3;q^5)_{\infty}}, \qquad |q| < 1.\label{rr2}
\end{align}
Rogers \cite{rogers-1894} showed that
\begin{equation}\label{rrcf-product}
R(q)=q^{1/5}\df{(q;q^5)_{\infty}(q^4;q^5)_{\infty}}{(q^2;q^5)_{\infty}(q^3;q^5)_{\infty}}=q^{1/5}\df{H(q)}{G(q)}, \qquad |q| < 1.
\end{equation}
Thus, $R(q)$ can be represented as a quotient of the Rogers--Ramanujan functions \eqref{rr1} and \eqref{rr2}.


The identities \eqref{rr1} and \eqref{rr2} were communicated by Ramanujan \cite[p.~344]{RamanujanCollected} in his first letter to Hardy on 16 January 1913. However, that page of Ramanujan's first letter has been lost! In a problem that Ramanujan \cite{1914b}, \cite[p.~330]{RamanujanCollected} submitted to the \emph{Journal of the Indian Mathematical Society} in 1914, he wrote, ``Examine the correctness of the following results,'' namely, \eqref{rr1} and \eqref{rr2}. Ramanujan's formulation appears to indicate that he did not have a proof of \eqref{rr1} and \eqref{rr2} at that time, and consequently also at the time that he wrote to Hardy. We conjecture that Hardy separated the page of Ramanujan's letter containing the identities \eqref{rr1} and \eqref{rr2} to show them to others, and consequently the separated page was never returned to its original position in Ramanujan's letter.

P.~A.~MacMahon stated the identities without proof in his famous book \cite{MM}, published in 1916, and devoted an entire chapter to them. Assuming their validity, he found combinatorial interpretations for both \eqref{rr1} and \eqref{rr2}. The first identity implies that the number of partitions of a positive integer $n$ into distinct nonconsecutive parts is equinumerous with the number of partitions of $n$ into parts that are congruent to either 1 or 4 modulo 5. The second identity implies that the number of partitions of a positive integer $n$ into distinct nonconsecutive parts, with each part at least 2, is equinumerous with the number of partitions of $n$ into parts that are congruent to either 2 or 3 modulo 5.

In 1917, while perusing old volumes of the \emph{Proceedings of the London Mathematical Society}, by chance, Ramanujan came across Rogers' paper \cite{rogers-1894} from 1894 proving the aforementioned identities. Ramanujan then found his own proof, and Rogers greatly simplified his former proof. Hardy arranged for both proofs to be published together \cite{RR1} in 1919. Since 1917, these identities have inspired an enormous amount of ongoing research involving other identities in the theory of $q$-series and their combinatorial interpretations.

Later, Hardy wrote \cite[p.~91]{12}:
``They [the Rogers--Ramanujan identities] were found first in 1894 by Rogers, a mathematician of great talent but comparatively little reputation, now remembered mainly from Ramanujan's rediscovery of his work. Rogers was a fine analyst, whose gifts were, on a smaller scale, not unlike Ramanujan's; but no one paid much attention to anything he did, and the particular paper in which he proved the formulae was quite neglected.''

In his first letter to Hardy, Ramanujan \cite[p.~xxviii]{RamanujanCollected}, \cite[p.~57]{bcbrar} tantalizingly informed him that he has established theorems on further continued fractions.

``The above theorem is a particular case of a theorem on the continued fraction
\begin{equation*}
\df{1}{1}\+\df{ax}{1}\+\df{ax^2}{1}\+\df{ax^3}{1}\+\df{ax^4}{1}\+\df{ax^5}{1}\+\cds,
\end{equation*}
which is a particular case of the continued fraction
\begin{equation*}
\df{1}{1}\+\df{ax}{1+bx}\+\df{ax^2}{1+bx^2}\+\df{ax^3}{1+bx^3}\+\cds,
\end{equation*}
which is a particular case of a general theorem on continued fractions.'' However, Ramanujan gave no information about these more general theorems.

In his second notebook \cite[p.~196]{nb}, \cite[p.~30]{III}, Ramanujan provided the information that he withheld from Hardy. Entry~15 of Chapter~16 in his second notebook is given by
\begin{equation}\label{entry15a}
\df{\displaystyle{\sum_{n=0}^{\infty}}\df{b^nq^{n^2}}{(aq;q)_n(q;q)_n}}{\displaystyle{\sum_{n=0}^{\infty}}\df{b^nq^{n(n+1)}}{(aq;q)_n(q;q)_n}} =
1+\df{bq}{1-aq}\+\df{bq^2}{1-aq^2}\+\df{bq^3}{1-aq^3}\+\cds , \qquad |q| < 1.
\end{equation}
In a corollary, Ramanujan records the special case of \eqref{entry15a} when $a = 0$, namely,
\begin{equation}\label{corollary,entry15a}
\df{\displaystyle{\sum_{n=0}^{\infty}}\df{b^nq^{n^2}}{(q;q)_n}}{\displaystyle{\sum_{n=0}^{\infty}}\df{b^nq^{n(n+1)}}{(q;q)_n}} = 1+\df{bq}{1}\+\df{bq^2}{1}\+\df{bq^3}{1}\+\cds, \qquad |q| < 1.
\end{equation}
If $b = 1$, \eqref{corollary,entry15a} reduces to the formula \eqref{rrcf-product} for the Rogers--Ramanujan continued fraction, upon use of the Rogers--Ramanujan identities \eqref{rr1} and \eqref{rr2}.

\section{Interlude: Ramanujan's Notebooks}

In the previous paragraph, we briefly mentioned Ramanujan's notebooks. Ramanujan who was born on December 22, 1887 in the southern Indian town, Erode, probably began to record his discoveries in a notebook around 1903 or 1904 when he entered the Government College in his home town of Kumbakonam. He worked on a slate, as paper was expensive. Thus, it was prohibitive for him to record his proofs in his notebook. After arriving in England in 1914, he concentrated on publishing his results, and from a letter that he wrote to a friend in Madras, at least at the beginning of his stay, Ramanujan apparently added few entries to his notebook. However, by the time that he returned to India in March, 1919, he had prepared a second notebook, which was a revised edition of his first notebook, but with many more new theorems. The second notebook contains 21 chapters totalling 256 pages, followed by 100 pages of unorganized results. Ramanujan also recorded his findings in a short third notebook of 33 pages. More information about the history and content of Ramanujan's notebooks \cite{nb} can be found in the first author's first book \cite{I} on Ramanujan's notebooks.

\section{Values of the Rogers--Ramanujan Continued Fraction}\label{values}

Recall from \eqref{S} that $S(q) = - R(-q)$. As illustrated by \eqref{2} and \eqref{3}, generally, if one can evaluate $R(e^{-2\pi\sqrt{n}})$, then one can evaluate $S(e^{-\pi\sqrt{n}})$ as well, and conversely. Moreover, from \eqref{symmetry}, observe that if we can determine the value of $R(e^{-2\alpha})$, then we can also find the value of $R(e^{-2\pi^2/\alpha})$.

In addition to the elegant evaluations found in \eqref{2}, \eqref{3}, and \eqref{5}, Ramanujan determined $R(e^{-\pi\sqrt{n}})$ for several additional values of $n$, where $n$ is a positive rational number. On page~311 in his first notebook \cite{nb}, \cite[p.~20]{V}, Ramanujan, with an ingenious argument, determined values for $R(q)$, when $q = e^{-4\pi}, e^{-8\pi}, e^{-16\pi}, e^{-6\pi}$, in identical formats. Define
\begin{equation}\label{6}
2c := 1+\df{a+b}{a-b}\sqrt{5},
\end{equation}
where $a$ and $b$ are distinct real numbers. If $a = 5^{1/4}$, $b = 1$, and $c$ is given by \eqref{6}, then
\begin{equation}\label{7}
R(e^{-4\pi}) = \sqrt{c^2+1} - c.
\end{equation}
If $a = (60)^{1/4}$, $b = 2-\sqrt{3}+\sqrt{5}$, and $c$ is given by \eqref{6}, then
\begin{equation}\label{8}
R(e^{-6\pi}) = \sqrt{c^2 + 1} - c.
\end{equation}
The values of $R(q)$ for $q = e^{-8\pi}$ and $q = e^{-16\pi}$ are similarly obtained \cite[pp.~20--30]{V}.
The first proof of \eqref{7} is due to K.~G.~Ramanathan \cite{ramanathan}. Using Ramanujan's parametrization \eqref{6}, Heng Huat Chan and the first author \cite{hhc1} gave the first proof for the set of all four values. They also determined some new values. For example \cite[p.~899]{hhc1}, \cite[p.~65]{geabcbI},
\begin{equation*}
S(e^{-\pi/\sqrt{3/5}}) = \left(\df{-5\sqrt{5}-3+\sqrt{30(5 + \sqrt{5})}}{4}\right)^{1/5}.
\end{equation*}
The first author, Chan, and L.-C.~Zhang \cite{hhc-lcz} extended Ramanujan's ideas in \eqref{6}. 

The authors of \cite{hhc1} and \cite{hhc-lcz}, and also Chan \cite{hhc-2}, found further new values of $R(q)$, which do not arise from the format \eqref{6}. In particular, they used \emph{modular equations} \cite[pp.~4--5, 213--214]{III} for several evaluations. For example, Chan \cite[pp.~344, 351--352]{hhc-2}, \cite[p.~156]{hhc-lcz} proved that
\begin{equation*}
S(e^{-\pi\sqrt{3}}) = \df{-3-\sqrt{5}+\sqrt{6(5 + \sqrt{5})}}{4}.
\end{equation*}

Next, we define a \emph{class invariant}. After Ramanujan \cite[p.~37]{III}, first set
\begin{equation*}
\chi(q) := (-q;q^2)_{\infty}, \qquad |q|<1.
\end{equation*}
When $q = e^{-\pi\sqrt{n}}$, for a positive rational number $n$, the class invariant $G_n$ is defined by \cite[pp.~21, 183]{V}
\begin{equation}\label{def:G}
G_n := 2^{-1/4}q^{-1/24}\chi(q).
\end{equation}
Class invariants play key roles in determining values of theta functions and continued fractions. Class invariants are algebraic numbers, and, generally, they are difficult to determine. These sources \cite{RamanujanModularPi}, \cite[pp.~23--39]{RamanujanCollected}, \cite[Chapter~34]{V} contain tables of values for $G_n$.

The Rogers--Ramanujan continued fraction lives in the world of theta functions. Because of the brevity of this survey, we shall not define a general theta function \cite[p.~34]{III}, but we will define one particular theta function. After Ramanujan \cite[p.~36, Entry~22]{III}, define the most well-known theta function $\varphi(q)$ by
\begin{equation}\label{phi}
\varphi(q) := \sum_{n=-\infty}^{\infty}q^{n^2}, \qquad |q|<1.
\end{equation}

In \cite[Theorem~5.2]{cubic-quintic}, the authors proved the following theorem for a class of quotients of theta functions \eqref{phi}. They proved a similar theorem with `5' replaced by `3' on the left side of \eqref{11}. See also \cite[pp.~330, 339]{V}.

\begin{theorem}\label{thm:3n} If $n$ is a positive rational number, then
\begin{equation}\label{11}
\frac{\varphi(e^{-5\pi\sqrt{n}})}{\varphi(e^{-\pi\sqrt{n}})} = \frac{1}{\sqrt{5}}\bigg(1 + \frac{2 G_{25n}}{G_{n}^5}\bigg)^{1/2}.
\end{equation}
\end{theorem}
Thus, if we know the values of the requisite class invariants on the right-hand side above, we can determine the value of the quotient of theta functions on the left-hand side.

We provide one example. From \cite[pp.~189--190]{V}, we know that
\begin{equation}\label{10}
G_1 = 1 \qquad \text{and} \qquad G_{25} = \df{\sqrt{5} + 1}{2}.
\end{equation}
Thus, setting $n = 1$ in \eqref{11} and employing \eqref{10}, after algebraic simplification, we find that \cite[p.~327]{V}
\begin{equation*}
\frac{\varphi(e^{-5\pi})}{\varphi(e^{-\pi})} = \frac{1}{(5\sqrt{5} - 10)^{1/2}}.
\end{equation*}

Following our paper \cite[Theorem~4.2]{cubic-quintic}, for $|q| < 1$, let
\begin{equation}\label{def:RR-p}
p := \frac{4q(-q; q^2)_\infty}{(-q^5; q^{10})_\infty^5} = \frac{4q\chi(q)}{\chi^5(q^5)}.
\end{equation}
The present authors \cite{rrcf1} established the following general theorem, which can be engaged to determine values of $R(q)$, provided that the requisite class invariants are known.

\begin{theorem}\label{thm:RR-quintic}
For $0 < |q|<1$,
\begin{equation*}
R(q) = \df{u}{\sqrt{p + 1} + 1} = \df{\sqrt{p + 1} - 1}{v} \qquad \text{and} \qquad 
R(q^4) = \df{\sqrt{p + 1} - 1}{u} = \df{v}{\sqrt{p + 1} + 1},
\end{equation*}
with
\begin{equation*}
u = (\alpha p)^{1/5} \qquad \text{and} \qquad v = (\beta p)^{1/5},
\end{equation*}
 where $\alpha$ and $\beta$ are roots of the polynomial equation
\begin{equation}\label{eq:RR-polynomial}
x^2 - ((p - 1)^2 + 7) x + p^2 = 0.
\end{equation}
\end{theorem}

The values of $u$ and $v$ depend on the order of $\alpha$ and $\beta$. This ambiguity can be resolved by using the definitions in \cite[Theorem~4.2(0)]{cubic-quintic}, but we omit these details for simplicity.

If $q = e^{-\pi\sqrt{n}}$, where $n$ is a positive rational number, then the value of $p$ in \eqref{def:RR-p} can be expressed in terms of the class invariants \eqref{def:G}, and the solutions of \eqref{eq:RR-polynomial}, in fact, lead to the values
\begin{equation*}
u,v = \left\{\df{p}{2}\left((p - 1)^2 + 7 \pm (4 - p)(4 + p^2)^{1/2}\right)\right\}^{1/5},
\end{equation*}
where $u$ is taken with the $+$ sign, and $v$ is taken with the $-$ sign.

Theorem~\ref{thm:RR-quintic} is an analogue of theorems in \cite[Entry~1, Theorem~2.1, Theorem~4.2]{cubic-quintic} and \cite[Entry~1.1, Theorem~4.1]{rebak} for calculating values of $\varphi(e^{-\pi\sqrt{n}})$. We can make use of Theorem~\ref{thm:RR-quintic} to calculate values of $R(e^{-\pi\sqrt{n}})$, when $n$ is a positive rational number. For example, if we take $q = e^{-\pi}$ in Theorem~\ref{thm:RR-quintic}, with the aid of \eqref{def:G} and \eqref{10}, 
we are led to a more explicit version of \eqref{7}, namely,
\begin{equation*}
R(e^{-4\pi}) = \frac{\sqrt{5} + 1}{2} \cdot \frac{(5\sqrt{5} - 10)^{1/2} - 1}{1 + (5 - 2\sqrt{5})^{1/2}}.
\end{equation*}

As a corollary of Theorem~\ref{thm:RR-quintic}, we note that
\begin{equation*}
\df{1}{R(q)}-R(q^4) = \df{2}{u} \qquad \text{and} \qquad \df{1}{R(q^4)}-R(q) = \df{2}{v}.
\end{equation*}

Chan \cite{hhc-2}, and the first author, Chan, and Zhang \cite{hhc-lcz} used class invariants to also determine specific values of 
$R(e^{-\pi\sqrt{n}})$. Further evaluations of $R(e^{-\pi\sqrt{n}})$ are derived in the authors' paper \cite{rrcf1}. See papers by Jinhee Yi \cite{yi-2} and by N.~D.~Baruah and N.~Saikia \cite{baruah-1} for additional evaluations.

\section{Identities Involving $R(q)$}

In his notebooks \cite{nb}, lost notebook \cite{lnb}, and previously unpublished manuscripts that are contained in \cite{lnb}, Ramanujan recorded a plethora of beautiful identities involving $R(q)$, $R(q^n)$, $n \in \mathbb{N}$, and theta functions at different arguments. We provide only a short sampling. Since the initial proofs of the identities from Ramanujan's notebooks and lost notebook, many of them have been given different proofs in the literature. Moreover, finding new identities continues to attract many researchers. Major sources for Ramanujan's identities involving $R(q)$ are Chapters~1--4 of the first book by Andrews and the first author on Ramanujan's lost notebook \cite{geabcbI}, Chapter~32 of the first author's book \cite{V}, and the Memoir \cite{memoir} by Andrews, the first author, L.~Jacobsen, and R.~L.~Lamphere.

One of the most important identities for $R(q)$ is given by
\begin{equation}\label{eq:identity-R-1}
\df{1}{R(q)}-1-R(q) = \df{(q^{1/5};q^{1/5})_{\infty}}{q^{1/5}(q^5;q^5)_{\infty}},
\end{equation}
which is found in Ramanujan's second notebook \cite[Chapter~16]{nb}, \cite[pp.~265--267]{III} and his lost notebook \cite{lnb}, \cite[pp.~11--18]{geabcbI}, and first proved by Watson \cite{watson1929a}. Ramanujan also proved a companion to \eqref{eq:identity-R-1}, namely \cite[p.~11]{geabcbI},
\begin{equation}\label{eq:identity-R-2}
\df{1}{R^5(q)}-11-R^5(q) = \df{(q;q)_{\infty}^6}{q(q^5;q^5)_{\infty}^6}.
\end{equation}
Soon-Yi Kang \cite{kang-acta} gave proofs of both \eqref{eq:identity-R-1} and \eqref{eq:identity-R-2}, and several related theorems as well, and then used her results to determine many values of $R(e^{-\pi\sqrt{n}})$ and $S(e^{-\pi\sqrt{n}})$, where $n$ is a positive rational number. Chadwick Gugg \cite{gugg} also gave an elegant proof of \eqref{eq:identity-R-2}.
Remarkably, Ramanujan \cite[p.~207]{lnb}, \cite[p.~13]{geabcbI} also saw that, in fact, his two identities above are special cases of very general two-variable identities involving theta functions!

On page~206 in his lost notebook \cite{lnb}, \cite[pp.~21--24]{geabcbI}, in addition to a set of very unusual theta function identities proved by Seung Hwan Son \cite{son}, \cite[pp.~180--194]{geabcbII} and the second author~\cite{rebak}, Ramanujan offers four insightful formulas that factorize \eqref{eq:identity-R-1} and \eqref{eq:identity-R-2}. We provide the first two of them, which arise from the factorization of \eqref{eq:identity-R-1}. In Ramanujan's notation, first let
\begin{equation*}
\alpha := \df{1 - \sqrt{5}}{2} \qquad \text{and} \qquad \beta :=\df{1 + \sqrt{5}}{2}.
\end{equation*}
Then, if $t=R(q)$,
\begin{align}
\df{1}{\sqrt{t}}-\alpha\sqrt{t} &= \df{1}{q^{1/10}}\sqrt{\df{(q;q)_{\infty}}{(q^5;q^{5})_{\infty}}}\prod_{n=1}^{\infty}\df{1}{1+\alpha q^{n/5}+q^{2n/5}},\label{eq:factorization-1}\\
\df{1}{\sqrt{t}}-\beta\sqrt{t} &= \df{1}{q^{1/10}}\sqrt{\df{(q;q)_{\infty}}{(q^5;q^{5})_{\infty}}}\prod_{n=1}^{\infty}\df{1}{1+\beta q^{n/5}+q^{2n/5}}.\label{eq:factorization-2}
\end{align}
Note that if we invert the roles of $\alpha$ and $\beta$ in \eqref{eq:factorization-1}, we obtain \eqref{eq:factorization-2}. Multiplying \eqref{eq:factorization-1} and \eqref{eq:factorization-2}, we deduce \eqref{eq:identity-R-1}.

On page~321 in his second notebook \cite{nb}, \cite[p.~17]{V} and on page~365 in his lost notebook \cite{lnb}, \cite[p.~92]{geabcbI}, \cite[p.~27]{memoir}, Ramanujan recorded a beautiful identity connecting $R(q)$ with $R(q^3)$. Let $u=R(q)$ and $v=R(q^3)$. Then,
\begin{equation*}
(v - u^3)(1 + uv^3) = 3u^2v^2.
\end{equation*}
He also stated \cite[p.~92]{geabcbI} similar identities connecting the pairs: $R(q), S(q)$; $R(q), R(q^4)$; and $R(q), R(q^5)$. Proofs for all of these identities have been provided by J.~Yi \cite{yi-1}. See also Gugg's paper \cite{gugg}.

In his second notebook \cite[p.~326]{nb}, Ramanujan introduced the parameter
\begin{equation*}
k := R(q)R^2(q^2)
\end{equation*}
and declared that
\begin{equation*}
R^5(q) = k\left(\df{1-k}{1+k}\right)^2 \qquad \text{and} \qquad R^5(q^2) = k^2\left(\df{1+k}{1-k}\right).
\end{equation*}
See \cite[Entry~24]{memoir} or \cite[pp.~12--14]{V} for proofs. In his lost notebook \cite[pp.~53, 56, 208]{lnb}, \cite[pp.~33--39, 81--84]{geabcbI}, Ramanujan continued his study of $k$. For example, if $k \leq \sqrt{5} -2$, then
\begin{equation*}
R(q^{1/2}) = \df{k^{1/10}(1 + k)^{4/5}(1 - k)^{1/5}}{\sqrt{k} + \sqrt{1 + k - k^2}}.
\end{equation*}
Except for one identity proved in \cite{hhc-lcz}, S.-Y.~Kang proved all of Ramanujan identities involving $k$ in her paper \cite{kang-ram}.

\section{Interlude: Ramanujan's Lost Notebook}

Ramanujan died on 26 April 1920. His wife, Janaki, subsequently donated all of her late husband's papers to the University of Madras. On 30 August 1923, Francis Dewsbury, Registrar at the University of Madras, sent to Hardy a packet of Ramanujan's papers. Unfortunately, neither Hardy nor Dewsbury made an account of what was sent. Sometime in the late 1930s, Hardy passed on this material to G.~N.~Watson, who earlier had devoted at least 10 years to editing Ramanujan's notebooks. Watson was Professor of Mathematics at the University of Birmingham. He died in 1965, and shortly thereafter, Robert Rankin, who had earlier succeeded Watson as Professor of Mathematics at the University of Birmingham, but who was then Professor of Mathematics at the University of Glasgow, paid a visit to Watson's widow. Watson had left piles of jumbled, unorganized papers, both mathematical and otherwise, on his attic floor office. In view of her late husband's fondness for Cambridge, Rankin and Mrs.~Watson agreed that during the next few years Rankin would sort through Watson's papers, and those of possible value were to be sent to Trinity College Library, Cambridge. Rankin found a considerable quantity of Ramanujan's manuscripts, most of which were incomplete or fragmentary, and sent them to Cambridge on 26 December 1968. Evidently, other than the librarians at Cambridge, Rankin did not inform anyone about his shipment.

George Andrews visited Cambridge in the spring of 1976. Cambridge mathematician, Lucy Slater, advised Andrews that he might find it interesting to browse through Watson's papers at Trinity College Library. In doing so, he found a pile of 138 sides of papers in Ramanujan's handwriting. Andrews immediately (and naturally) named this clump of papers, \emph{Ramanujan's Lost Notebook}. Subsequently, other papers were located in the libraries at both Cambridge and Oxford. All of these were gathered together and published in one volume \cite{lnb} in 1988.

\section{Finite Form of the Rogers--Ramanujan Continued Fraction}

In the Introduction, we wrote ``Infinite continued fractions are considerably more interesting than finite continued fractions.'' Perhaps that claim was too hastily given. There are finite forms of the Rogers--Ramanujan continued fraction given in Entry~16 of Chapter~16 in Ramanujan's second notebook \cite[p.~196]{nb}, \cite[p.~31]{III}.

For each positive integer $n$, let
\begin{equation*}
\mu := \mu_n(a, q) := \sum_{k=0}^{\lfloor (n+1)/2 \rfloor}\df{a^kq^{k^2}(q;q)_{n-k+1}}{(q;q)_k(q;q)_{n-2k+1}}
\end{equation*}
and
\begin{equation*}
\nu := \nu_n(a, q) := \sum_{k=0}^{\lfloor n/2 \rfloor}\df{a^kq^{k(k+1)}(q;q)_{n-k}}{(q;q)_k(q;q)_{n-2k}}.
\end{equation*}
Then,
\begin{equation}\label{finite-rr}
\df{\mu}{\nu} = 1 + \df{aq}{1}\+\df{aq^2}{1}\+\cds\+\df{aq^n}{1}.
\end{equation}
Note the similarity of \eqref{finite-rr} with \eqref{rrcf-product}, and also with \eqref{corollary,entry15a}.

\section{An Uncommon Continued Fraction `Identity'}

We return to \eqref{cf2}, which can be loosely construed as an analogue of the Rogers--Ramanujan continued fraction, wherein the powers of $q$ are replaced by the powers themselves. In his second letter to Hardy \cite[p.~58]{bcbrar}, and on page~276 in his second notebook \cite{nb}, \cite[p.~67]{V}, Ramanujan offered an unusual `identity.'
``When $x$ is small,
\begin{multline}
\df{1}{1}\+\df{1}{1}\+\df{2}{1}\+\df{3}{1}\+\df{4}{1}\+\cds = x\sqrt{e}\sum_{n=1}^{\infty}e^{-(1+nx)^2/2}\\
+\df{x}{2}-\df{x^2}{12}-\df{x^4}{360}-\df{x^6}{5040}-\df{x^8}{60480}-\df{x^{10}}{1710720} \qquad \text{nearly.''}\label{interesting}
\end{multline}
Consider the product on the right-hand side of \eqref{interesting}. Ramanujan wanted Hardy to examine it ``when $x$ is small.'' As $x \to 0^+$, the terms in the infinite series approach $e^{-1/2}$, and so the infinite series tends to $\infty$, as $x \to 0$, or non-rigorously, if we set $x = 0$ in the summands, we obtain a divergent series. Also, on the right-hand side is a factor $x\sqrt{e}$, which, of course, tends to $0$ as $x \to 0^+$. Thus, on the right side, as $x$ tends to $0^+$, one factor tends to $0$, while the other tends to $\infty$. Which of them is going to win as $x \to 0^+$ in the limit? Well, it's a tie. The product does not tend to $0$ or to $\infty$, but instead it tends to the infinite continued fraction on the left-hand side of \eqref{interesting}. But Ramanujan tells us more. When $x > 0$ is small, he gives a better approximation. If we move the polynomial in $x$ to the left-hand side, for $x > 0$ very, very small, the polynomial on the left-hand side is negative. Thus, the now modified left-hand side is slightly less than the value of the continued fraction. Hence, as $x \to 0^+$, the infinite series is `ahead' of $x\sqrt{e}$, but in the end they tie. In his notebooks, Ramanujan offers many approximations, especially in powers of the variable $x$, and he often appended his approximations with the words, ``nearly,'' or ``very nearly.''

The wonderful formula \eqref{interesting} is, in fact, a special case of a more general result in Ramanujan's second notebook \cite[p.~276]{nb}, \cite[p.~4]{memoir}, \cite[p.~67]{V}. In his generalization, Ramanujan gives a formula for the coefficients in the associated polynomial. If we let the degree of the polynomial tend to infinity, we obtain an infinite series that does \emph{not} converge. Watson \cite{watson1931} proved a version of Ramanujan's claim in 1931.

The continued fraction \eqref{cf2} has no known elementary evaluation \cite[p.~423]{finch}. However, Ramanujan \cite{1914a}, \cite[p.~329]{RamanujanCollected} posed the following remarkable problem in the \emph{Journal of the Indian Mathematical Society} in 1914:
\begin{equation*}
\bigg(1 + \frac{1}{1 \cdot 3} + \frac{1}{1 \cdot 3 \cdot 5} + \frac{1}{1 \cdot 3 \cdot 5 \cdot 7} + \cdots\bigg) +
\bigg(\df{1}{1}\+\df{1}{1}\+\df{2}{1}\+\df{3}{1}\+\df{4}{1}\+\cds\bigg) = \sqrt{\frac{\pi e}{2}}.
\end{equation*}
A more general identity is in Ramanujan's second notebook \cite[p.~152]{nb}, \cite[p.~166, Entry~43]{II}.

\end{document}